\input amstex
\documentstyle{amsppt}
\define\tr{\text{trace}}
\define\fp{(\text{f.p.})_{z=0}}

\define\reso{(\lambda-a)^{-1}}
\define\res{\text{res}}
\NoRunningHeads
\topmatter
\title On multiplicative properties of determinants\endtitle
\author Leonid Friedlander \endauthor
\affil University of Arizona\endaffil
\endtopmatter
\document
\head 1. Introduction\endhead
Let $A$ be an elliptic classical pseudodifferential operator of positive order $k$ on a closed compact manifold $M$ of
dimension $d$. Suppose that $A$ admits an Agmon angle, which is a solid angle $\{\lambda\in\Bbb C: \alpha\leq \arg\lambda\leq\beta\}$ that is free
from values of the principal symbol of $A$. We  assume that the null space of $A$ is trivial. Then one can define complex powers $A^z$ of $A$ and the logarithm $\log A$.
The operator $\log A$ is not a classical PDO but an operator with $\log$-polyhomogeneous symbol ([L]); actually, the only
logarithmic term in the expansion of its symbol is $\log |\xi|$. It is well known that the function
$$\zeta(z)=\tr A^z,\tag 1.1$$
wich is a priori defined and analytic in the half-plane $\text{Re} z,-d/k$, admits an analytic continuation to a meromorphic
function in the whole complex plane; the point $z=0$ is a regular point. Usually, the zeta-function is taken as function of 
$s=-z$. Then the determinant of $A$ is defined according to the formula
$$\log\det A=\zeta'(0).\tag 1.2$$

In the half-plane ${\text Re}z<-d/k$ one can differentiate (1): $\zeta'(z)=\tr (\log A)A^z$, so
$$\log\det A=\tr (\log A)A^z\big|_{z=0}.\tag 1.3$$
The expression on the right in (1.3) can be interpreted as a regularized trace of $\log A$; the operator $A$ is used 
as the regularizer. It is known that the determinant of elliptic operators is not multiplicative; however the multiplicative
anomaly $\log\det (AB)-\log\det A-\log\det B$ can be computed in terms of symbols of the operators $A$ and $B$ ([F], [KV]).
In this paper, I discuss a different situation: what happens if one takes an operator $I+T$ for $B$ where $T$ is an operator
of negative order. This problem may arise when one wants to compare determinants of two elliptic operators
with the same principal symbol. In the case when the order of $T$ is smaller than $-d$, the operator $T$ is of trace class,
and it is not difficult to show that $\det(A(I+T))=\det A\det(I+T)$ (e.g., see [F]). In general, the operator $T$ 
belongs to the Schatten class ${\frak S}_p$ for $p>d/s$ if $T$ is a PDO of order $-s$.

Let us recall the notion of a regularized Fredholm determinant (e.g., see [GK], [S]). Let $T$ be a compact operator in
a Hilbert space, and $T\in{\frak S}_p$. Let $m\geq p$ be an integer. Then
$$\text{det}_m(I+T)=\prod_{j=1}^\infty\biggl((1+\lambda_j)\exp \sum_{p=1}^{m-1}\frac{(-1)^p}{p}\lambda_j^p\biggr)\tag 1.4$$
where $\lambda_j$ are eigenvalues of $T$; each eigenvalue is counted as many times as its multiplicity is.
The product on the right in (1.4) converges. One can re-write (1.4) in the form
$$\log\text{det}_m(I+T)=\tr\biggl(\log(I+T)+\sum_{p=1}^{m-1}\frac{(-1)^p}{p}T^p\biggr).\tag 1.5$$

From this point, we assume that $T$ is a classical pseudodifferential operator and the operator
$I+T$ is invertible.
If the operator $T$ is not of trace class, there is no reason to expect the determinant of $A(I+T)$ to be equal to
$\det A\hbox{\rm det}_m(I+T)$. My goal is to find an expression for 
$$w_m(A,T)=\log\det(A(I+T))-\log\det A-\log\hbox{\rm det}_m(I+T).\tag 1.6$$

\head 2. Variational formulas\endhead
Let $w(t)=w_m(A,tT)$. Clearly, $w(0)=0$ and $w(1)=w_m(A,T)$. The variational formula for the determinant of an elliptic 
operator is well known:
$$\frac{d}{dt}\log\det (A(t))=\fp\tr(A'(t)A(t)^{z-1})\tag 2.1$$
where $A'(t)$ is the derivative of the operator-valued function $A(t)$.
The meromorphic function $\tr(A'(t)A(t)^{z-1})$ is either regular at $z=0$ or $z=0$ is its simple pole; the expression
on the right in (2.1) is the constant term of its Laurent expansion at $z=0$. Therefore,
$$\align\frac{d}{dt}\log\det(A(I+tT))&=\fp\tr[AT(A(I+tT))^{z-1}]\\
                                     &=\fp\tr[AT(A(I+tT))^z(I+tT)^{-1}A^{-1}]\\
                                     &=\fp\tr[(I+tT)^{-1}T(A(I+tT))^z]\\
                                     &=\sum_{p=1}^{m-1}(-1)^{p-1}t^{p-1}\fp\tr[T^p(A(I+tT))^z]\\&+
                                     (-1)^{m-1}t^{m-1}\tr [T^m(I+tT)^{-1}].\tag 2.2
\endalign$$
Here we assume the operator $T^m$ to be of trace class. Formula (1.5) implies
$$\align
\frac{d}{dt}\log\text{\rm det}_m(I+tT)&=\tr[T(I+tT)^{-1}-\sum_{p=1}^{m-1}(-1)^{p-1}t^{p-1}T^p]\\
 &=(-1)^{m-1}t^{m-1}\tr [T^m(I+tT)^{-1}].\tag 2.3
 \endalign$$
From (2.2) and (2.3) we conclude
$$w'(t)=\sum_{p=1}^{m-1}(-1)^{p-1}t^{p-1}\fp\tr[T^p(A(I+tT))^z].\tag 2.4$$

The operator valued function
$$\Phi_p(t,z)=T^p\frac{(A(I+tT))^z-A^z}{z}$$
is a holomorphic operator valued function in the sense of Guillemin [G]; operators $\Phi_p(t,z)$ are classical PDO of order $kz-ps$ (recall that $T$ is an operator of order $-s$.) The function $\tr \Phi_p(t,z)$ that is initially defined in the half-plane $\text{Re}(z)<(-d+ps)/k$ admits an analytic continuation to a meromorphic function in the whole complex plane, and it has simple poles. In particular, $z=0$ may be a pole. Then,
$$\res_{z=0}\tr\Phi_p(t,z)=\frac{1}{k}\res Q_p(t)\tag 2.5$$
(see [G], [W]). Here $\res Q_p(t)$ is the Guillemin-Wodzicki non-commutative residue of
the operator $Q_p(t)=\Phi_p(t,0)$. The non-commutative residue of an operator $Q$
$$\res Q=\int_{S^*M}q_{-d}(x,\xi)\mu\tag 2.6$$
where $\mu=\alpha \wedge (d\alpha)^{d-1}$ with $\alpha=\xi dx$ being the Liouville form on $T^*M$, and $q_{-d}(x,\xi)$  is the homogeneous of degree $-d$ term in the complete symbol expansion of the operator $Q$. This term is not invariantly defined as a function on $T^*M$, so one has to use local coordinates and a partition of unity to get an expression in the right in (2.5). A remarkable fact that is due to Guillemin and Wodzicki is that the result is independent of the choice made. It follows from (2.5) that
$$\fp\tr[T^p(A(I+tT))^z]=\fp\tr(T^pA^z)+\frac{1}{k}\res Q_p(t),$$
and formula (2.4) can be re-written as
$$w'(t)=\sum_{p=1}^{m-1}(-1)^{p-1}t^{p-1}\fp\tr(T^pA^z)+\frac{1}{k}\sum_{p=1}^{m-1}(-1)^{p-1}t^{p-1}\res Q_p(t).\tag 2.7$$
The first sum on the right in (2.7) is a polynomial in  $t$ and hence can be integrated explicitly; the second sum contains local expressions only.

Notice that 
$$Q_p(t)=T^p[\log(A(I+tT))-\log A].\tag 2.8$$
Operators $\log (A(I+tT))$ and $\log A$ are not classical PDO; they are PDO with poly-logarithmic symbols; however their difference is a classical PDO. The order of $Q_p$ is
$-(p+1)s$. Keeping in mind that $w_m(A,T)=w(1)$ and $w(0)=0$ one gets
\proclaim{Theorem}
Let $A$ be an elliptic pseudo-differential operator of positive order $k$ on a $d$-dimensional closed manifold that posesses an Agmon angle. Let $T$ ba a classical pseudo-differential operator on $M$ 
of negative order $-s$. 
Assume that the null space of $A$ is trivial and the operator $I+T$ is invertible. Let $m$ be an integer number that is greater than $d/s$. Then
$$\align \log\det (A(I+T))-\log\det A-\log\text{det}_m(I+T) &
=\sum_{p=1}^{m-1}\frac{(-1)^{p-1}}{p}\fp\tr(T^pA^z)\\
&+\frac{1}{k}\sum_{p=1}^{m-1}(-1)^{p-1}\int_0^1 t^{p-1}\res Q_p(t)dt \tag 2.9
\endalign$$
where operators $Q_p$ are given in (2.8) and $\res$ is the Guillemin--Wodzicki non-commutative residue.
\endproclaim
\head 3. The structure of $\res Q_p(t)$ and some special cases\endhead
Let
$$a(x,\xi)\sim \sum_{j=0}^\infty a_{k-j}(x,\xi)\quad\text{and}\quad
\tau(x,\xi)\sim\sum_{j=0}^\infty\tau_{-s-j}(x,\xi)$$
be complete symbol expansions for operators $A$ and $T$ in local coordinates. An index shows the degree of homogeneity of the corresponding function. By $(\lambda-a)^{-1}$ 
I denote the complete symbol of the parametrix of the operator $\lambda I-A$ as a pseudo-differential operator with parameter (e.g., see [Sh]).  It is the inverse of $\lambda-a$ 
in the non-commutative symbolic algebra. In the symbolic algebra of pseudo-differential operators with parameter, $\lambda$ is treated as an additional dual variable, and one assigns weight $k$ to it. Then
$$(\lambda-a(1+t\tau))^{-1}-(\lambda-a)^{-1}=
\sum_{q=1}^\infty t^q(\lambda-a)^{-1}[a\tau(\lambda-a)^{-1}]^q\tag 3.1$$
is the symbol of the difference between resolvents of the operators $A(I+tT)$ and $A$.
All multiplications on the right in (3.1) are symbolic algebra multiplications.
The $q$-th term in the sum in (3.1) is of order $-k-qs$, so the sum makes sense as an asymptotic sum. The principal symbol of the difference between resolvents equals
$a_k(x,\xi)\tau_{-s}(x,\xi)/(\lambda-a_k(x,\xi))^2$, and all homogeneous terms in the symbolic expansion of the symbol (3.1) are rational functions of $\lambda$ with the only pole at $\lambda=a_k(x,\xi)$; 
the residue at this pole vanishes.

Let $\ell$ be a ray in the complex plane that goes from $0$ to $\infty$ and lies in the Agmon angle for $A$. We make a cut in the complex plane along $\ell$ to define $\log\lambda$. Let $\Gamma$ be the 
contour that 
goes from $-\infty$ to the point $\epsilon$-close to $0$
along one side of the cut, then goes around $0$ to the opposite side of the cut, and then goes back to $\infty$; here $\epsilon<|a_k(x,\xi)|$. Then
$$\log(a(1+t\tau))-\log a=\sum_{q=1}^\infty t^q\frac{1}{2\pi i}
\int_\Gamma\log\lambda\reso [a\tau\reso]^qd\lambda.\tag 3.2$$
Once more, all operations in (3.2) are being done in the symbolic algebra sense. Despite of the presence of $\log\lambda$ in the integrands in (3.2) the result is a classical symbol because in each 
homogeneous term the residue at the only pole $\lambda=a_k(x,\xi)$ vanishes. The order of the $q$-th term in (3.2) equals $-qs$.Finally, for the symbop $q^{(p)}(t,x,\xi)$ of the operator $Q_p(t)$ one gets
$$q^{(p)}(t, x,\xi)=\sum_{j=1}^\infty t^j\tau^j\frac{1}{2\pi i}
\int_\Gamma\log\lambda\reso [a\tau\reso]^jd\lambda.\tag 3.3$$
The order of the $q$-th term on the right in (3.3) equals $-(q+p)s$; its residue vanishes if $q+p>d/s$ because the residue of a symbol of order smaller than $-d$ equals $0$.
Therefore,
$$\res Q_p(t)=\sum_{j=1}^{[(d/s)-p]} t^j\res\biggl(\tau^p\frac{1}{2\pi i}
\int_\Gamma\log\lambda\reso [a\tau\reso]^jd\lambda.\biggr)\tag 3.4$$
Here, by $[\alpha]$ I denote the floor of the number $\alpha$.

Let us discuss two special cases. The first one is when $s>d/2$. An example would be the difference $\log\det(\Delta+u(x))-\log\det\Delta$ on manifolds of dimension $2$ or $3$.
Though the operator $\Delta$ is not invertible, its 
(modified) determinant equals the determinant of $A=\Delta+P$ where $P$ is the orthogonal projection onto the null space of $\Delta$. One can set $T=(u-P)(\Delta+P)^{-1}$.  Notice that $P$ is a pseudo-differential operator of order $-\infty$, and, for the purpose of computing symbols, can be disregarded. In this case $m=2$, there is just one operator, $Q_1$ on the right in (2.9), its order equals $-2s<-d$, so its residue vanishes. One gets
$$\det(A(I+T))-\det A-\hbox{det}_2(I+T)=\fp\tr (TA^z).$$

The second case is when $s=d/2$. Then $m=3$. The operator $Q_2$ is of order $-3s<-d$,
so its residue vanishes. However, the operator $Q_1$ is of order $-2s=-d$, and
$$\res Q_1(t)=\int_{S^*(M)}q^{(1)}_{-d}(t, x,\xi)\mu.$$
Only the $q=1$ term in (3.4) contributes to the principle symbol of $Q_1$. By evaluating the corresponding integral one gets
$$q^{(1)}_{-d}(t, x,\xi)=t\tau_{-d/2}(x,\xi)^2.$$
Then formuls (2.9) takes the form
$$\align\det(A(I+T))&-\det A-\hbox{det}_3(I+T)=\fp\tr (TA^z)\\ 
&-\frac{1}{2}\fp\tr(T^2A^z)
+\frac{1}{2k}\int_{S^*M}\tau_{-d/2}(x,\xi)^2\mu.
\endalign$$


\Refs
\widestnumber\key{KVa}
\ref \key F \by L. Friedlander
\book Determinants of elliptic operators
\bookinfo Ph.D. Dissertation
\publ MIT\publaddr Cambridge, MA\yr 1989
\endref
\ref\key FG \by L.Friedlander, V. Guillemin
\paper Determinants of zeroth order operators
\jour J. Diff. Geometry\yr 2008\vol 78\pages 1--12
\endref
\ref \key G 
\paper A new proof of Weyl's formula on the asymptotic distribution of eigenvalues
\by V. Guillemin
\jour Advances in Math.\yr 1985\vol 55\pages 131--160
\endref
\ref\key GK \by I. Gohberg, M. G. Krein
\book Introduction into the theory of non-selfadjoint operators
\publ Springer\publaddr Berlin-Heidelberg-New York
\yr 1971
\endref
\ref \key KV\by M. Kontsevich, S. Vishik
\book Determinants of elliptic pseudo=differential operators
\bookinfo preprint\publ Max-Plank Institut f\"ur Mathematik
\publaddr Bonn\yr 1994
\endref
\ref \key L\by M. Lesch
\paper On the non-commutative residue for pseudodifferential operators with \newline  log-polyhomogeneous symbols
\jour Ann. Global Anal. Geom.\yr 1999\vol 17\pages 151--187
\endref
\ref\key S\by B. Simon
\book Trace Ideals and their Applications, Second Edition
\bookinfo  Mathematical Surveys and Monographs
\publ American Mathematical Society
\publaddr Providence
\yr 2010
\endref
\ref \key W\by M. Wodzicki
\book Noncommutative residue\bookinfo Lect. Notes in Math. 1289
\publ Springer-Verlag\yr 1987
\endref




\endRefs
\enddocument